\documentclass{amsart}
\usepackage{amsmath,amssymb,amsthm,mathrsfs}
\usepackage{enumitem}
\usepackage[colorlinks,linkcolor=blue,citecolor=blue]{hyperref}
\allowdisplaybreaks
\numberwithin{equation}{section}

\newtheorem{theorem}{Theorem}[section]
\newtheorem{proposition}[theorem]{Proposition}
\newtheorem{corollary}[theorem]{Corollary}
\newtheorem{lemma}[theorem]{Lemma}
\theoremstyle{definition}

\newtheorem{remark}[theorem]{Remark}

\newcommand{\ep}{\varepsilon}
\newcommand{\Om}{\Omega}
\newcommand{\la}{\lambda}
\newcommand{\R}{\mathbb{R}}
\newcommand{\II}{\mathrm{II}}

\begin{document}

\title[A sharp threshold in the Li--Lin open problem]{A sharp threshold in the Li--Lin open problem}

	\author{Qihan He}

\author{Zongyan Lv}

\address{College of Mathematics \& Center for Applied Mathematics of Guangxi (Guangxi University), Guangxi University, Nanning 530003, P. R. China}
\email{heqihan277@gxu.edu.cn}

\address{School of Mathematics and Physics, University of Science and Technology Beijing, Beijing 100083, P. R. China}
\email{zongyanlv0535@163.com}

\begin{abstract}
	In this paper, we study positive solutions to the problem
	\begin{equation*}
		-\Delta u = |x|^{-s_2}u^{2^*(s_2)-1} - \la |x|^{-s_1}u^{p-1} \quad \text{in }\Om, \qquad u=0 \quad \text{on }\partial\Om,
	\end{equation*}
	where $\Om\subset\R^N$ ($N\ge 3$) is a bounded domain with $0\in\partial\Om$, $0\le s_1<s_2<2$, $2^*(s)=\frac{2(N-s)}{N-2}$, $2^*(s_2)<p<2^*(s_1)$, and $\la>0$. This is the Li--Lin open problem \cite{LL} in the form studied by Cerami, Zhong, and Zou \cite{CZZ} and by Tang and Tang \cite{TT1}. Assume that $\partial\Om$ is $C^{1,\alpha}$ near $0$ with a leading term of homogeneity $1+\alpha$, and set $\theta=\frac{(N-2)(2^*(s_1)-p)}{2}$. We prove that a solution exists for every $\la>0$ if $\theta>\alpha$ and the $\alpha$-mean curvature of $\partial\Om$ at $0$ is negative, whereas $\la^*=\sup\{\la>0:\text{a solution exists}\}$ is finite if $\theta\le\alpha$, without any curvature assumption. For $\alpha=1$, this removes the restriction to almost every $\la$ in \cite[Theorem 1.6]{CZZ}, answers a question of Zou \cite{Z}, and disproves the conjecture in \cite[Remark 1.7]{TT1}. The proofs rely on a scale-invariant estimate uniform in $\la$, a Hadamard-type identity for the half-space extremal and a localized Pohozaev identity combined with sharp decay estimates. Notably, no regularity of $\partial\Om$ is required away from $0$.
\end{abstract}


\maketitle

\section{Introduction and main results}

Let $N\ge3$, $\Om\subset\R^N$ be a bounded domain with $0\in\partial\Om$ and 
\[
0\le s_1<s_2<2,\qquad q:=2^*(s_2)=\frac{2(N-s_2)}{N-2}<p<2^*(s_1).
\]
We consider the following equation
\begin{equation}\label{eq:P}
-\Delta u=|x|^{-s_2}u^{q-1}-\la|x|^{-s_1}u^{p-1},\qquad u>0\ \text{in}\ \Om,\qquad u\in H^1_0(\Om),\tag{$P_\la$}
\end{equation}
with $\la>0$, and set
\[
\la^*:=\sup\{\la>0:\eqref{eq:P}\ \text{has a solution}\}\in[0,\infty].
\]

In \cite{LL}, Li and Lin  studied elliptic equations with two Hardy--Sobolev terms and asked for the existence of positive solutions when $s_1<s_2$ and the term with the larger exponent has a negative sign. Cerami, Zhong and Zou \cite{CZZ} gave the first answers in the regime of \eqref{eq:P}: if $\Om$ is a $C^1$ domain whose boundary is $C^2$ at $0$ with mean curvature $H(0)<0$, then \eqref{eq:P} has a solution for small $\la$ \cite[Theorem 1.5]{CZZ}, and for almost every $\la>0$ provided $p<2^*(s_1)-\frac2{N-2}$ \cite[Theorem 1.6]{CZZ}. Zou \cite{Z} asked whether ``almost every'' can be removed and what happens for larger $p$. Tang and Tang \cite{TT1} proved that $\la^*<\infty$ in the doubly critical case $p=2^*(s_1)$, obtained very weak solutions for every $\la\in(0,\la^*)$ when $p$ is close to $2^*(s_1)$, and conjectured \cite[Remark 1.7]{TT1} that $\la^*=\infty$ whenever $p<2^*(s_1)$. For domains which are starshaped with respect to $0$, \eqref{eq:P} has no solution at all \cite{TT2}.

Throughout the paper, we set
\[
a:=\frac{N-2}{2},\ \ \Phi_u(x):=|x|^au(x),\ \ \theta:=\frac{(N-2)(2^*(s_1)-p)}{2}\in(0,s_2-s_1),\ \ c_p:=\frac{\theta}{p}.
\]
No regularity of $\partial\Om$ away from $0$ is used. Near $0$ we assume the following:
\begin{itemize}[leftmargin=3.2em]
\item[$(H_\alpha)$] There exist $r_\Om>0$, $\alpha\in(0,1]$ and $h\in C^{1,\alpha}(\R^{N-1})$ with $h(0)=0$, $\nabla h(0)=0$, such that, after a rotation, $\Om\cap B_{2r_\Om}=\{x_N>h(x')\}\cap B_{2r_\Om}$.
\item[$(D_\alpha)$] There is $h_{1+\alpha}\in C^1(\R^{N-1}\setminus\{0\})$, homogeneous of degree $1+\alpha$, such that $h=h_{1+\alpha}+o(|x'|^{1+\alpha})$ and $\nabla h=\nabla h_{1+\alpha}+o(|x'|^{\alpha})$ as $x'\to0$.
\end{itemize}
Under $(D_\alpha)$ we define the $\alpha$-mean curvature of $\partial\Om$ at $0$ by
\[
H_\alpha(0):=\frac{1}{|S^{N-2}|}\int_{S^{N-2}}h_{1+\alpha}(\sigma)\,d\sigma .
\]
If $\partial\Om$ is of class $C^2$ near $0$, then $(D_1)$ holds with $h_2(x')=\frac12\II_0(x',x')$ and $H_1(0)=\frac12 H(0)$, where $H(0)$ denotes the mean curvature of $\partial\Om$ at $0$ (with the sign convention of \cite{CZZ}).

Since $\theta>\alpha$ if and only if
\[
p<p_\sharp(\alpha):=2^*(s_1)-\frac{2\alpha}{N-2},
\]
our main result reads as follows.

\begin{theorem}\label{thm:main}
Assume $(H_\alpha)$. Then the following hold.
\begin{enumerate}[label=(\roman*)]
\item If $p<p_\sharp(\alpha)$, $(D_\alpha)$ holds and $H_\alpha(0)<0$, then \eqref{eq:P} has a solution for every $\la>0$, that is, $\la^*=\infty$.
\item If $p\ge p_\sharp(\alpha)$, then $\la^*<\infty$, without any assumption on the curvature of $\partial\Om$.
\end{enumerate}
In particular, if $\partial\Om$ is $C^{1,1}$ near $0$, twice differentiable at $0$ and $H(0)<0$, then $p_\sharp(1)=2^*(s_1)-\frac{2}{N-2}$ is a sharp threshold: $\la^*=\infty$ for $p<p_\sharp(1)$, and $\la^*<\infty$ for $p\ge p_\sharp(1)$.
\end{theorem}

For $\alpha=1$, part (i) removes ``almost every'' in \cite[Theorem 1.6]{CZZ}, and part (ii) answers the second question of \cite{Z}; together with \cite[Theorem 1.5]{CZZ}, it gives $0<\la^*<\infty$ for $p\ge p_\sharp(1)$ under the assumptions of \cite{CZZ}. The range of $p$ in \cite[Theorem 1.4]{TT1} satisfies $\theta\le\frac{s_2-s_1}{2}<1$, so part (ii) shows that the conjecture \cite[Remark 1.7]{TT1} fails throughout that range. For $\alpha<1$ the threshold $p_\sharp(\alpha)$ depends on the order of contact of $\partial\Om$ with its tangent plane at $0$.

The key ingredient in both parts is an a priori estimate for $\Phi_u$ which is uniform in $\la\ge\la_->0$ and requires no regularity of $\partial\Om$ away from $0$ (Proposition \ref{prop:univ}). We also describe the behaviour at $\la^*$, see Proposition \ref{prop:star}. If $\theta<\alpha$, then any sequence of solutions with $\la_n\to\la^*$ converges to a nontrivial $u_0$, and $\la^*$ is attained if and only if there are such solutions with bounded energy.

Let us comment on the proofs. The two nonlinear terms balance at the profile $U^*_\la=\la^{-1/(p-q)}|x|^{-\alpha_*}$, $\alpha_*=\frac{s_2-s_1}{p-q}$, and a blow-up argument alone fails since rescalings may converge to the constant solution of the limiting equation. We exclude such plateaus by an interior dichotomy and a connectedness argument, and then exclude the remaining profiles by Liouville theorems. In  $\R^n$,  we use the theorem of Li and Zhang \cite{LZ}. For Theorem \ref{thm:main}(i), the first order term in the energy expansion is computed by the identity of Lemma \ref{lem:hadamard}, which is a Hadamard formula for the half-space extremal, and Proposition \ref{prop:univ} supplies the $H^1$ bounds needed to pass from almost every $\la$ to every $\la$. For Theorem \ref{thm:main}(ii), a Pohozaev identity at a fixed radius, combined with a barrier adapted to $C^{1,\alpha}$ boundaries, gives $\la\le C\rho^{\alpha-\theta}$, where $\rho$ is the radius of the region in which $\Phi_u$ is not small.

\section{Preliminaries}\label{sec:pre}

We write
\[
\alpha_*=\frac{s_2-s_1}{p-q},\quad L_\la(x)=\la|x|^\theta,\quad U^*_\la=\la^{-\frac{1}{p-q}}|x|^{-\alpha_*},\quad \Phi^*=|x|^aU^*_\la=L_\la^{-\frac{1}{p-q}},
\]
and $\ell^*_\la=|x|(\Phi^*)^{-\frac{q-2}{2}}$. We shall use $a(q-2)=2-s_2$, $aq=N-s_2$ and
\[
\theta=N-s_1-ap=(\alpha_*-a)(p-q)=2-s_1-a(p-2).
\]
We also set $d(x)=\operatorname{dist}(x,\partial\Om)$, $c_0=a^{2/(q-2)}$ and
\begin{align*}
	A(u,E)=\int_{\Om\cap E}|x|^{-s_1}u^p, \quad A(u)=A(u,\R^N), \quad  B(u)=\int_\Om|x|^{-s_2}u^q,\\
I_\la(u)=\frac12\int_\Om|\nabla u|^2+\frac{\la}{p}A(u_+)-\frac1qB(u_+),\qquad G(x,u)=\frac{u^q}{q|x|^{s_2}}-\frac{\la u^p}{p|x|^{s_1}},\\
\text{and}\quad P(r,u)=\int_{\Om\cap\partial B_r}\Big[r|\partial_ru|^2-\frac r2|\nabla u|^2+rG(x,u)+au\,\partial_ru\Big]d\sigma .
\end{align*}
Under $(H_\alpha)$, decreasing $r_\Om$, we may assume that
\begin{align*}
	|\nabla h(x')|\le K|x'|^\alpha, \quad |h(x')|\le K|x'|^{1+\alpha}\quad \text{and} \quad K(2r_\Om)^\alpha\le\frac14.
\end{align*}
  On $\partial\Om\cap B_{2r_\Om}$ the outer normal $\nu$ satisfies
\begin{equation}\label{eq:xnu}
x\cdot\nu=\frac{x'\cdot\nabla h(x')-h(x')}{\sqrt{1+|\nabla h(x')|^2}},\qquad |x\cdot\nu|\le C|x|^{1+\alpha},
\end{equation}
and, under $(D_\alpha)$, Euler's identity $x'\cdot\nabla h_{1+\alpha}=(1+\alpha)h_{1+\alpha}$ gives $x\cdot\nu=\alpha h_{1+\alpha}(x')+o(|x'|^{1+\alpha})$.

\begin{lemma}\label{lem:reg}
Let $u$ be a solution of \eqref{eq:P}. Then the following hold.
\begin{enumerate}[label=(\roman*)]
\item $u\in L^\infty(\Om)\cap C^2(\Om)$.
\item Under $(H_\alpha)$, $u\in C^{1,\gamma}(\overline\Om\cap B_{r_\Om})$ for some $\gamma>0$, and $u\le Cd$ in $\Om\cap B_{r_\Om}$.
\item Under $(H_\alpha)$, for $0<r\le r_\Om$,
\begin{equation}\label{eq:poh}
P(r,u)=-\la c_pA(u,B_r)-\frac12\int_{\partial\Om\cap B_r}(x\cdot\nu)|\partial_\nu u|^2\,d\sigma .
\end{equation}
\end{enumerate}
\end{lemma}

\begin{proof}
(i) Since $-\Delta u\le|x|^{-s_2}u^{q-1}$, the Brezis--Kato iteration \cite{BK}, with test functions $u\min\{u,L\}^{2k}$ and the Hardy--Sobolev inequality, gives $u\in L^t(\Om)$ for all $t<\infty$. Then $f:=|x|^{-s_2}u^{q-1}\in L^t$ for some $t>\frac N2$, and $w:=\Gamma*(f1_\Om)$ is bounded. Since $0\le(u-w)_+\le u$, we have $(u-w)_+\in H^1_0(\Om)$, and $-\Delta(u-w)\le0$ gives $u\le w$.

(ii) Rounding the edge of $\Om\cap B_{7r_\Om/4}$ we obtain a $C^{1,\alpha}$ domain $D$ with $\Om\cap B_{3r_\Om/2}\subset D\subset\Om\cap B_{2r_\Om}$. Let $\tilde h$ be harmonic in $D$ with boundary values $\|u\|_\infty\psi$, where $\psi\in C^\infty$, $\psi=0$ in $B_{5r_\Om/4}$ and $\psi=1$ outside $B_{4r_\Om/3}$. Then $u\le G_D[f]+\tilde h$ and $\tilde h\le Cd$ in $B_{r_\Om}$. Since $G_D(x,y)\le Cd(x)|x-y|^{1-N}$ in \cite{GW}, the bound $u\le C|y|^{\gamma_k}$ implies $u\le Cd$ if $t_k:=s_2-(q-1)\gamma_k<1$, and $u\le Cd|x|^{1-t_k}$ if $1<t_k<2$. Starting from $\gamma_0=0$ and iterating $\gamma_{k+1}=2-s_2+(q-1)\gamma_k$ (slightly decreased if $t_k=1$), we reach $t_k<1$ after finitely many steps. Then $|\Delta u|\le C|x|^{q-1-s_2}$, and \cite[Theorem 8.33]{GT} gives the $C^{1,\gamma}$ regularity.

(iii) Using $\frac{N-s_2}{q}=a$ and $\frac{N-s_1}{p}-a=c_p$, we see that
\[
\operatorname{div}\Big[(x\cdot\nabla u)\nabla u-\frac12|\nabla u|^2x+Gx+au\nabla u\Big]=-\la c_p|x|^{-s_1}u^p\quad\text{in}\ \Om .
\]
We integrate over $\{x_N>h_k(x')+\frac1k\}\cap(B_r\setminus\overline B_\delta)$, where $h_k\to h$ in $C^1$ are smooth, and let $k\to\infty$ and $\delta\to0$. On $\partial\Om$ we have $\nabla u=(\partial_\nu u)\nu$ and $G=0$, and by (ii) the integral over $\partial B_\delta$ is $O(\delta^N)$.
\end{proof}

\begin{lemma}\label{lem:KO}
For every bounded domain $\Om$, every $\la>0$ and every solution of \eqref{eq:P}, we have
\[
\Phi_u\le C_1\big(L_\la^{-\frac{1}{p-q}}+L_\la^{-\frac{1}{p-2}}\big),
\]
with $C_1=C_1(N,s_1,s_2,p)$. In particular, $u\le2C_1U^*_\la$ wherever $\Phi^*\ge1$.
\end{lemma}

\begin{proof}
We first note that if $w\in H^1_0(\Om)$, $w\ge0$ and $\Delta w\ge g$ in $\mathcal D'(\Om)$ with $g\in L^1_{\rm loc}(\R^N)$, then the zero extension $\tilde w$ satisfies $\Delta\tilde w\ge g1_{\{w>0\}}$ in $\mathcal D'(\R^N)$. Indeed, for $0\le\varphi\in C^\infty_c(\R^N)$ and $T_\ep=\min\{w/\ep,1\}$ we have $\varphi T_\ep\in H^1_0(\Om)$ and
\[
\int T_\ep\nabla w\cdot\nabla\varphi=\int\nabla w\cdot\nabla(\varphi T_\ep)-\frac1\ep\int_{\{w<\ep\}}\varphi|\nabla w|^2\le-\int g\varphi T_\ep ;
\]
it remains to let $\ep\to0$.

Fix $x_0\in\Om$, $B=B(x_0,|x_0|/2)$ and $ m=2^{\frac{1}{p-q}+\alpha_*}\la^{-\frac{1}{p-q}}|x_0|^{-\alpha_*}.$
On $\{u>m\}\cap B$ we have $\la|x|^{-s_1}u^{p-1}\ge2|x|^{-s_2}u^{q-1}$, and hence $w:=(u-m)_+\in H^1_0(\Om)$ satisfies 
\begin{align*}
	\Delta\tilde w\ge c\la|x_0|^{-s_1}\tilde w^{p-1} \quad \text{in} ~~ B.
\end{align*}
 Let $R=|x_0|/2$, $\sigma=\frac{2}{p-2}$ and $W=A_0(R^2-|x-x_0|^2)^{-\sigma}$. Since $\sigma(p-1)=\sigma+2$, a direct computation shows that
 \begin{align*}
 	\Delta W\le c\la|x_0|^{-s_1}W^{p-1} \quad \text{in} \quad  B \qquad \text{for}\quad  A_0^{p-2}=C(N,p)R^2/(\la|x_0|^{-s_1}).
 \end{align*}
   Since $\tilde w$ is bounded, $(\tilde w-W)_+\in H^1_0(B)$, and the maximum principle gives $\tilde w\le W$. Thus $w(x_0)\le W(x_0)\le C(\la|x_0|^{2-s_1})^{-\frac{1}{p-2}}$, and multiplying by $|x_0|^a$ we obtain the claim. If $\Phi^*\ge1$, then $L_\la^{-1/(p-2)}=(\Phi^*)^{\frac{p-q}{p-2}}\le\Phi^*$.
\end{proof}

\begin{lemma}\label{lem:nondeg}
Every solution of \eqref{eq:P} satisfies $\sup\limits_\Om\Phi_u\ge c_0$.
\end{lemma}

\begin{proof}
By the Nehari identity and the Hardy inequality,
\[
\|\nabla u\|_2^2\le B(u)=\int_\Om|x|^{-2}u^2\Phi_u^{q-2}\le(\sup\Phi_u)^{q-2}a^{-2}\|\nabla u\|_2^2.\qedhere
\]
\end{proof}

\begin{lemma}\label{lem:scaled}
Assume $(H_\alpha)$, $|x_0|\le r_\Om$ and $\Phi_u\le M$ in $\Om\cap B(x_0,|x_0|/2)$. Then $-\Delta u=Vu$ there with $|V|\le C(M)|x_0|^{-2}$, and
\[
\|\nabla u\|_{L^\infty(\Om\cap B(x_0,|x_0|/4))}\le C(M,\Om)|x_0|^{-1}\sup_{\Om\cap B(x_0,|x_0|/2)}u .
\]
\end{lemma}

\begin{proof}
We have $V=|x|^{-2}(\Phi_u^{q-2}-L_\la\Phi_u^{p-2})$. If $L_\la\le1$, then $L_\la\Phi_u^{p-2}\le M^{p-2}$. If $L_\la\ge1$, Lemma \ref{lem:KO} gives $L_\la\Phi_u^{p-2}\le CL_\la^{1-\frac{p-2}{p-q}}\le C$. After scaling by $|x_0|$ the boundary has bounded $C^{1,\alpha}$ norm, and \cite[Theorem 8.33]{GT} applies.
\end{proof}

\section{A uniform estimate}\label{sec:univ}

\begin{proposition}\label{prop:univ}
Assume $(H_\alpha)$ and let $\la_->0$. Then there exists $C_0>0$, depending only on $N,s_1,s_2,p,\Om$ and $\la_-$, such that every solution of \eqref{eq:P} with $\la\ge\la_-$ satisfies $\Phi_u\le C_0$ in $\Om$.
\end{proposition}

For $\kappa\ge0$ let $f_\kappa(t)=t^{q-1}-\kappa t^{p-1}$ and for $\kappa>0$ let $z_\kappa=\kappa^{-1/(p-q)}$.

\begin{lemma}\label{lem:liouRn}
Let $n\ge2$, and if $n\ge3$ assume $p<\frac{2n}{n-2}$. If $v\ge0$ is a bounded $C^2$ solution of $-\Delta v=f_\kappa(v)$ in $\R^n$, then either $v\equiv0$, or $\kappa>0$ and $v\equiv z_\kappa$.
\end{lemma}

\begin{proof}
Let $v\not\equiv0$; then $v>0$. If $\kappa>0$, translating along a sequence $v(x_k)\to\sup v$ we see that $f_\kappa(\sup v)\ge0$, that is, $v\le z_\kappa$. If $n=2$, then $v$ is a bounded positive superharmonic function, hence a constant $c$ with $f_\kappa(c)=0$. Let $n\ge3$ and $\sigma=\frac{n+2}{n-2}$. Set $F=f_\kappa$ on $(0,z_\kappa]$ and $F=0$ on $(z_\kappa,\infty)$ (for $\kappa=0$, $F(t)=t^{q-1}$). For $0<t<z_\kappa$, since $\kappa t^{p-q}\le1$ and $q<p<\sigma+1$,
\[
\frac{d}{dt}\big(F(t)t^{-\sigma}\big)=t^{q-2-\sigma}\big[(q-1-\sigma)+(\sigma+1-p)\kappa t^{p-q}\big]\le(q-p)t^{q-2-\sigma}<0 .
\]
Thus $F(t)t^{-\sigma}$ is nonincreasing on $(0,\infty)$ and is not constant on the range of $v$. By the Liouville theorem of Li and Zhang \cite[Theorem 1.3]{LZ}, $v$ is a constant $c$ with $F(c)=0$, and hence $\kappa>0$ and $c=z_\kappa$.
\end{proof}

Since $p<2^*(s_1)\le2^*<\frac{2(N-1)}{N-3}$, Lemma \ref{lem:liouRn} applies in $\R^N$ and in $\R^{N-1}$.

\begin{lemma}\label{lem:liouhalf}
Let $v\ge0$ be a bounded solution of $-\Delta v=f_\kappa(v)$ in $\R^N_+$, $v\in C^2(\R^N_+)\cap C(\overline{\R^N_+})$, $v=0$ on $\partial\R^N_+$. Then either $v\equiv0$, or $\kappa>0$, $\partial_Nv>0$ and $v(x',t)\uparrow z_\kappa$ as $t\to\infty$.
\end{lemma}

\begin{proof}
If $v\not\equiv0$, then $v>0$, and the monotonicity theorem for bounded positive solutions in half-spaces with $f(0)\ge0$ (see \cite{BCN2,BCN,Da}) gives $\partial_Nv>0$. Then $\lim_{t\to\infty}v(\cdot,t)$ is a bounded positive solution in $\R^{N-1}$, and Lemma \ref{lem:liouRn} applies.
\end{proof}

Let $g(x)=x_N-h(x')$. Since $\operatorname{Lip}(h)\le1$, we have $g/\sqrt2\le d\le g$ in $\Om\cap B_{r_\Om}$. Moreover, if $r'\le r_\Om$ and $0\le\tau(t)\le t/4$, then the set
\[
D=\{x\in\Om:|x|<r',\ g(x)>\tau(|x|)\}
\]
is connected and contains the segment $\{te_N:0<t<r'\}$. Indeed, along the arc $y(\phi)=|x|(\sin\phi\,\omega,\cos\phi)$ joining $x$ to $|x|e_N$, we have
\[
\frac{d}{d\phi}g(y(\phi))\le-|x|\sin\phi+K|x|^{1+\alpha}\sin^\alpha\phi,
\]
which is negative unless $\sin\phi\le s_*:=(2K|x|^\alpha)^{1/(1-\alpha)}$ (with $s_*=0$ if $\alpha=1$), while $g\ge|x|/2>\tau(|x|)$ when $\sin\phi\le s_*$.

\begin{proof}[Proof of Proposition \ref{prop:univ}]
By Lemma \ref{lem:KO}, $\Phi_u\le C(r)$ on $\{|x|\ge r\}$ for $\la\ge\la_-$, so it suffices to work near $0$. Let $\delta=\frac18$.

\emph{Step 1.} There are $K_1,K_2'>0$ such that, for every $\la>0$, every solution and every $x$ with $\Phi^*(x)\ge K_2'$ and $d(x)\ge K_1\ell^*(x)$,
\begin{equation}\label{eq:dich}
u(x)/U^*(x)\in[0,\delta)\cup(1-\delta,1+\delta).
\end{equation}
Otherwise there are $(\la_n,u_n,x_n)$ violating \eqref{eq:dich} with $K_{1,n},K_{2,n}'\to\infty$. Let $v_n(z)=u_n(x_n+\ell^*_nz)/U^*(x_n)$. Since $\ell^*_n/|x_n|=\Phi^*(x_n)^{-(q-2)/2}\to0$, the function $v_n$ is defined in $B_{K_{1,n}}$, it satisfies $v_n\le2^{1+\alpha_*}C_1$ by Lemma \ref{lem:KO}, and
\[
-\Delta v_n=\Big(\frac{|y|}{|x_n|}\Big)^{-s_2}v_n^{q-1}-\Big(\frac{|y|}{|x_n|}\Big)^{-s_1}v_n^{p-1},\qquad y=x_n+\ell^*_nz .
\]
A subsequence converges in $C^1_{\rm loc}(\R^N)$ to a bounded solution of $-\Delta v=f_1(v)$, so $v\in\{0,1\}$ by Lemma \ref{lem:liouRn}, which contradicts $v(0)\in[\delta,1-\delta]\cup[1+\delta,\infty)$.

\emph{Step 2.} Choose $K_2\ge\max\{K_2',(\la_-r_\Om^\theta)^{-1/(p-q)}\}$ with $2K_1K_2^{-(q-2)/2}\le\frac14$. Let 
\begin{align*}
	r_\la=(\la K_2^{p-q})^{-1/\theta} \quad \text{and} \quad  I_\la^0=\{x\in\Om:|x|<r_\la,\ g(x)>2K_1\ell^*(x)\}.
\end{align*}
Since $\ell^*$ is radial and $\Phi^*>K_2$ in $B_{r_\la}$, we have $2K_1\ell^*(x)\le|x|/4$ there, so $I_\la^0$ is connected, contains the axis, and $d>K_1\ell^*$ on it. By Step 1 and continuity, either $u<\delta U^*$ on $I^0_\la$, or $u>(1-\delta)U^*$ on $I^0_\la$. The latter forces $u\to\infty$ along the axis, which contradicts $u\in L^\infty(\Om)$. Hence $u<\delta U^*$ on $I^0_\la$.

\emph{Step 3.} Suppose that $\la_n\ge\la_-$ and $\Phi_{u_n}(y_n)\to\infty$; then $y_n\to0$. Let
\[
X=\overline\Om\cap\overline B_{r_\Om/2},\quad \Gamma=\{0\}\cup(\overline\Om\cap\partial B_{r_\Om/2}),\quad M_n=|x|^{-\frac{s_2}{2}}u_n^{\frac{q-2}{2}},\quad \tilde M_n=M_n+|x|^{-1}.
\]
By Lemma \ref{lem:reg}(ii), $\tilde M_n$ is positive and continuous in $X\setminus\Gamma$, and $M_n|x|=\Phi_{u_n}^{(q-2)/2}$. The doubling lemma \cite[Lemma 5.1]{PQS}, with $k_n=\frac14\tilde M_n(y_n)|y_n|\to\infty$, gives points $x_n'$ such that
\[
\tilde M_n(x_n')\operatorname{dist}(x_n',\Gamma)>2k_n\qquad\text{and}\qquad\tilde M_n\le2\tilde M_n(x_n')\ \ \text{in}\ X\cap B\big(x_n',k_n/\tilde M_n(x_n')\big).
\] Since $\operatorname{dist}(x_n',\Gamma)\le|x_n'|$, we get $\Phi_{u_n}(x_n')\to\infty$, and by Lemma \ref{lem:KO}, $L_{\la_n}(x_n')\to0$, so that $|x_n'|\to0$ and $\Phi^*(x_n')\to\infty$. Let $M_n'=M_n(x_n')$, so that $M_n\le3M_n'$ in $B(x_n',k_n/(2M_n'))\subset B(x_n',|x_n'|/2)$, and
\[
v_n(z)=\frac{u_n(x_n'+z/M_n')}{u_n(x_n')},\qquad -\Delta v_n=\Big(\frac{|y|}{|x_n'|}\Big)^{-s_2}v_n^{q-1}-\kappa_n\Big(\frac{|y|}{|x_n'|}\Big)^{-s_1}v_n^{p-1},
\]
where $\kappa_n=(u_n/U^*)^{p-q}(x_n')\le(2C_1)^{p-q}$. Then $v_n(0)=1$, $v_n$ is bounded on $B_{k_n/2}$, and $\ell^*(x_n')M_n'=\kappa_n^{\frac{q-2}{2(p-q)}}$. Up to a subsequence, $\kappa_n\to\kappa$, the domains $M_n'(\Om-x_n')$ converge to $\R^N$ or to a half-space $\{z\cdot\nu>-h_0\}$ with $h_0>0$, and $v_n\to v$ in $C^1_{\rm loc}$.

If the limit domain is $\R^N$, Lemma \ref{lem:liouRn} and $v(0)=1$ give $\kappa=1$. Since $d(x_n')M_n'\to\infty$, we have $x_n'\in I^0_{\la_n}$ for large $n$, while $u_n/U^*(x_n')=\kappa_n^{1/(p-q)}\to1$, which contradicts Step 2.

If the limit domain is a half-space, Lemma \ref{lem:liouhalf} gives $\kappa>0$ and $v(-h_0\nu+t\nu)\to z_\kappa$. Choose $t_0\ge4K_1\kappa^{\frac{q-2}{2(p-q)}}$ with $v(-h_0\nu+t_0\nu)\ge(1-\frac\delta2)z_\kappa$, and let $\tilde x_n=x_n'+(t_0-h_0)\nu_n/M_n'$, where $\nu_n$ is the inner normal at a boundary point nearest to $x_n'$. Then
\begin{gather*}
\frac{u_n}{U^*}(\tilde x_n)\to\kappa^{\frac{1}{p-q}}v(-h_0\nu+t_0\nu)\ge1-\frac\delta2,\quad 
d(\tilde x_n)M_n'\to t_0,\quad \ell^*(\tilde x_n)M_n'\to\kappa^{\frac{q-2}{2(p-q)}},
\end{gather*}
so that $\tilde x_n\in I^0_{\la_n}$ and $u_n>(1-\delta)U^*$ at $\tilde x_n$, which again contradicts Step 2.
\end{proof}

\begin{corollary}\label{cor:univ}
Assume $(H_\alpha)$ and $\la\ge\la_-$. Then $|\nabla u|\le C|x|^{-N/2}$ in $\Om\cap B_{r_\Om/2}$, and for $r\le r_\Om/2$, 
\[
A(u,B_r)\le Cr^\theta,\qquad\int_{\partial\Om\cap B_r}|x\cdot\nu||\partial_\nu u|^2\le Cr^\alpha,\qquad A(u)\le C_A .
\]
\end{corollary}

\begin{proof}
This follows from Proposition \ref{prop:univ}, Lemma \ref{lem:scaled} and \eqref{eq:xnu}, since $|x|^{-s_1}u^p\le C|x|^{\theta-N}$ and $|x\cdot\nu||\partial_\nu u|^2\le C|x|^{1+\alpha-N}$.
\end{proof}

\section{Existence for every \texorpdfstring{$\la>0$}{lambda>0}}\label{sec:exist}

Let
\[
\mu:=\inf\Big\{\|\nabla u\|_2^2:\ u\in D^{1,2}_0(\R^N_+),\ \int_{\R^N_+}|y|^{-s_2}|u|^q=1\Big\},\qquad\Theta:=\Big(\frac12-\frac1q\Big)\mu^{\frac{q}{q-2}}.
\] By \cite{Eg,GR1}, $\mu$ is attained by a positive function $w$ depending only on $(|y'|,y_N)$, with
\[
-\Delta w=\mu|y|^{-s_2}w^{q-1}\ \ \text{in}\ \R^N_+,\qquad w=0\ \ \text{on}\ \partial\R^N_+,\qquad\int_{\R^N_+}|y|^{-s_2}w^q=1 .
\]
We first record the regularity and decay of $w$.

\begin{lemma}\label{lem:w}
The minimizer $w$ satisfies
\begin{equation}\label{eq:wdecay}
w\in C^1(\overline{\R^N_+}),\qquad w(y)\le Cy_N(1+|y|)^{-N},\qquad|\nabla w(y)|\le C(1+|y|)^{-N}.
\end{equation}
\end{lemma}

\begin{proof}
\emph{Step 1.} We show that $w\in L^\infty(B_R^+)$ for every $R>0$. Write $-\Delta w=Vw$ with $V=\mu|y|^{-s_2}w^{q-2}$ and split $V=V1_{\{w>M\}}+V1_{\{w\le M\}}$. By the H\"older inequality with respect to $|y|^{-s_2}dy$ and the Hardy--Sobolev inequality, we obtain
\[
\int V1_{\{w>M\}}\varphi^2\le\mu\Big(\int_{\{w>M\}}|y|^{-s_2}w^q\Big)^{\frac{q-2}{q}}\Big(\int|y|^{-s_2}|\varphi|^q\Big)^{\frac2q}\le\delta(M)\|\nabla\varphi\|_2^2,
\]
with $\delta(M)\to0$ as $M\to\infty$. Since $s_2<2$, for every $\ep>0$ we have $|y|^{-s_2}\le\ep|y|^{-2}+C_\ep$, and the Hardy inequality gives
\[
\int V1_{\{w\le M\}}\varphi^2\le\mu M^{q-2}\Big(\frac{4\ep}{(N-2)^2}\|\nabla\varphi\|_2^2+C_\ep\|\varphi\|_2^2\Big).
\]
Hence the Brezis--Kato argument \cite{BK}, with test functions $\zeta^2w\min\{w,L\}^{2k}$ and a cut-off $\zeta$, gives $w\in L^t(B_R^+)$ for every $t<\infty$. Choosing $t>\frac N2$ with $s_2t<N$, we get $\mu|y|^{-s_2}w^{q-1}\in L^t(B_R^+)$, and the boundary local maximum principle \cite[Theorem 8.25]{GT} gives $w\in L^\infty(B_{R/2}^+)$.

\emph{Step 2.} The argument of Lemma \ref{lem:reg}(ii), applied in a smooth domain $D$ with $B_1^+\subset D\subset B_2^+$, gives $w\le Cy_N$ in $B_1^+$ and $w\in C^{1,\gamma}(\overline{B_1^+})$. Away from $0$ the coefficient $|y|^{-s_2}$ is smooth, so $w\in C^1(\overline{\R^N_+}\setminus\{0\})$.

\emph{Step 3.} Let $w^\star(y)=|y|^{2-N}w(y/|y|^2)$. The Kelvin transform preserves the Dirichlet integral and $\R^N_+$, and the change of variables $z=y/|y|^2$ gives
\[
\int|y|^{-s_2}|w^\star|^q\,dy=\int|z|^{s_2+(N-2)q-2N}|w|^q\,dz=\int|z|^{-s_2}|w|^q\,dz,
\]
since $(N-2)q=2N-2s_2$. Thus $w^\star$ is again a minimizer, and Steps 1 and 2 give $w^\star\le Cy_N$ in $B_1^+$. For $|y|\ge1$ it follows that
\[
w(y)=|y|^{2-N}w^\star(y/|y|^2)\le C|y|^{2-N}\frac{y_N}{|y|^2}=Cy_N|y|^{-N}.
\]
Together with Step 2, this gives $w\le Cy_N(1+|y|)^{-N}$. Finally, let $|y|=R\ge2$. In $B(y,R/2)\cap\R^N_+$ we have $-\Delta w=Vw$ with
\[
|V|\le C|y|^{-s_2}(R^{1-N})^{q-2}\le CR^{-2},
\]
since $s_2+(N-1)(q-2)\ge s_2+2(2-s_2)\ge2$. Rescaling $B(y,R/2)$ to the unit ball and using \cite[Theorem 8.33]{GT}, we obtain $|\nabla w(y)|\le CR^{-1}\sup_{B(y,R/2)}w\le CR^{-N}$.
\end{proof}

For $\beta>0$ we set $K_\beta(w)=\int_{\R^{N-1}}|y'|^{\beta}|\partial_Nw(y',0)|^2dy'$.

\begin{lemma}\label{lem:hadamard}
Let $g\in C^1(\R^{N-1}\setminus\{0\})$ be homogeneous of degree $1+\alpha$, $\alpha\in(0,1]$. Then
\begin{multline*}
q\int_{\R^N_+}\partial_Nw\,\nabla'w\cdot\nabla g\,dy-\mu s_2\int_{\R^N_+}|y|^{-s_2-2}y_Ng\,w^q\,dy\\
=-\frac q2\int_{\R^{N-1}}|\partial_Nw(y',0)|^2g(y')\,dy' .
\end{multline*}
\end{lemma}

\begin{proof}
By \eqref{eq:wdecay} all integrals converge absolutely, and the boundary terms on $\partial B_R$ below vanish as $R\to\infty$. We compute $\int_{\R^N_+}(-\Delta w)g\,\partial_Nw$ in two ways. Since $w=0$ on $\{y_N=0\}$ and $\partial_N|y|^{-s_2}=-s_2|y|^{-s_2-2}y_N$, we have
\[
\int(-\Delta w)g\,\partial_Nw=\frac\mu q\int|y|^{-s_2}g\,\partial_N(w^q)=\frac{\mu s_2}{q}\int|y|^{-s_2-2}y_Ng\,w^q .
\]
On the other hand, since $\partial_Ng=0$ and $\nabla w=(\partial_Nw)e_N$ on $\{y_N=0\}$, we have
\begin{align*}
\int\Delta w\,g\,\partial_Nw&=\int\operatorname{div}(g\,\partial_Nw\nabla w)-\frac12\int g\,\partial_N|\nabla w|^2-\int\partial_Nw\,\nabla'w\cdot\nabla g\\
&=-\frac12\int_{\R^{N-1}}|\partial_Nw|^2g\,dy'-\int\partial_Nw\,\nabla'w\cdot\nabla g .
\end{align*}
Comparing the two identities, we obtain the claim.
\end{proof}

Let $\eta\in C^\infty_c(B_{2r})$ with $\eta=1$ on $B_r$, $r\le r_\Om/2$, let $\varphi(x)=(x',x_N-h(x'))$, and for $0<\ep<1$ let
\[
u_\ep(x)=\eta(x)\ep^{-a}w\big(\varphi(x)/\ep\big)\in H^1_0(\Om).
\]

\begin{lemma}\label{lem:expansion}
Assume $(H_\alpha)$ and $(D_\alpha)$. Then, as $\ep\to0$,
\begin{gather*}
\|\nabla u_\ep\|_2^2=\mu-2\ep^\alpha J_1+o(\ep^\alpha),\qquad J_1=\int_{\R^N_+}\partial_Nw\,\nabla'w\cdot\nabla h_{1+\alpha}\,dy,\\
B(u_\ep)=1-s_2\ep^\alpha J_2+o(\ep^\alpha),\qquad J_2=\int_{\R^N_+}|y|^{-s_2-2}y_Nh_{1+\alpha}w^q\,dy,\\
A(u_\ep)=\ep^\theta\big(A(w)+o(1)\big),\qquad A(w)=\int_{\R^N_+}|y|^{-s_1}w^p\,dy .
\end{gather*}
Consequently,
\begin{equation}\label{eq:level}
\sup_{t\ge0}\Big[\frac{t^2}{2}\|\nabla u_\ep\|_2^2-\frac{t^q}{q}B(u_\ep)\Big]=\Theta+\frac{q\Theta}{(q-2)\mu}H_\alpha(0)K_{1+\alpha}(w)\,\ep^\alpha+o(\ep^\alpha).
\end{equation}
\end{lemma}

\begin{proof}
We change variables by $x'=\ep y'$, $x_N=\ep y_N+h(\ep y')$, whose Jacobian is $1$. Since
\begin{align*}
	 \nabla_x[w(\varphi/\ep)]=\ep^{-1}(\nabla'w-\partial_Nw\nabla h,\partial_Nw) \quad \text{and} \quad N-2a-2=0,
\end{align*}
 we have
\[
\|\nabla u_\ep\|_2^2=\int_{\R^N_+}\big[|\nabla w|^2-2\partial_Nw\,\nabla'w\cdot\nabla h(\ep y')+|\partial_Nw|^2|\nabla h(\ep y')|^2\big]dy+O(\ep^{N}),
\]
where the error comes from $\eta$ and \eqref{eq:wdecay}. By $(D_\alpha)$,
\begin{align*}
	\nabla h(\ep y')=\ep^\alpha\nabla h_{1+\alpha}(y')+o(\ep^\alpha)|y'|^\alpha \quad \text{and} \quad |\nabla h(\ep y')|\le K\ep^\alpha|y'|^\alpha, 
\end{align*}
 and dominated convergence gives the first expansion. 
Next, $|x|^2=\ep^2|y|^2(1+e)$, where
\[
e=\frac{2y_Nh(\ep y')/\ep+(h(\ep y')/\ep)^2}{|y|^2}=2\ep^\alpha\frac{y_Nh_{1+\alpha}(y')}{|y|^2}+o(\ep^\alpha)|y|^\alpha,\qquad|e|\le\frac12
\]
on the support of $\eta(\varphi^{-1}(\ep\cdot))$ for $r$ small. Since $N-aq=s_2$ and $N-ap=s_1+\theta$, we obtain
\[
B(u_\ep)=\int\eta^q|y|^{-s_2}(1+e)^{-\frac{s_2}{2}}w^q\,dy,\qquad A(u_\ep)=\ep^\theta\int\eta^p|y|^{-s_1}(1+e)^{-\frac{s_1}{2}}w^p\,dy,
\]
and dominated convergence gives the remaining expansions. Finally, the supremum in \eqref{eq:level} equals
\[
\Big(\frac12-\frac1q\Big)\|\nabla u_\ep\|_2^{\frac{2q}{q-2}}B(u_\ep)^{-\frac{2}{q-2}},
\]
which is
\[
\Theta-\frac{2\Theta}{(q-2)\mu}(qJ_1-\mu s_2J_2)\ep^\alpha+o(\ep^\alpha).
\]
By Lemma \ref{lem:hadamard} with $g=h_{1+\alpha}$, and since $|\partial_Nw(y',0)|^2$ is radial,
\[
qJ_1-\mu s_2J_2=-\frac q2\int_{\R^{N-1}}|\partial_Nw(y',0)|^2h_{1+\alpha}(y')\,dy'=-\frac q2H_\alpha(0)K_{1+\alpha}(w).\qedhere
\]
\end{proof}

For $\alpha=1$, \eqref{eq:level} is the expansion behind \cite[Lemma 2.4]{CZZ}. Lemma \ref{lem:hadamard} determines the sign of the first order term for every order of contact.

\begin{lemma}\label{lem:ps}
Assume $(H_\alpha)$, $(D_\alpha)$, $H_\alpha(0)<0$ and $\theta>\alpha$, and let $\la_0>0$. Then there exist $\bar c<\Theta$ and a set $\Lambda\subset[\la_0/2,2\la_0]$ of full measure such that, for every $\la\in\Lambda$, $I_\la$ has a bounded Palais--Smale sequence at a level $c_\la\in(0,\bar c]$.
\end{lemma}

\begin{proof}
We recall the theorem of Jeanjean \cite[Theorem 1.1]{J}. Let $X$ be a Banach space, $J\subset(0,\infty)$ an interval, and $\Upsilon_m(v)=\mathcal A(v)-m\mathcal B(v)$, $m\in J$, a family of $C^1$ functionals with $\mathcal B\ge0$ and $\mathcal A(v)\to\infty$ as $\|v\|\to\infty$. Assume that there are $v_1,v_2\in X$ such that, for every $m\in J$,
\[
\tilde c_m:=\inf_{\gamma\in\mathcal G}\max_{t\in[0,1]}\Upsilon_m(\gamma(t))>\max\{\Upsilon_m(v_1),\Upsilon_m(v_2)\},
\]
where $\mathcal G=\{\gamma\in C([0,1],X):\gamma(0)=v_1,\ \gamma(1)=v_2\}$.
Then, for almost every $m\in J$, there is a bounded sequence $(v_k)\subset X$ with $\Upsilon_m(v_k)\to\tilde c_m$ and $\Upsilon_m'(v_k)\to0$.

Let
\[
c_1=\frac{q\Theta}{(q-2)\mu}|H_\alpha(0)|K_{1+\alpha}(w)>0,\qquad T_0=(4q\mu)^{\frac{1}{q-2}},\qquad T=4^{\frac{1}{p-2}}T_0 .
\] By Lemma \ref{lem:expansion}, we notice that
\begin{align*}
	  \frac12\mu\le\|\nabla u_\ep\|_2^2\le2\mu, \quad  B(u_\ep)\ge\frac12 \quad \text{and} \quad A(u_\ep)\le2A(w)\ep^\theta \quad \text{for}  ~~\ep ~~ \text{small}.
\end{align*}
 Therefore, for every $\la\le2\la_0$ and every $s\in[0,T]$, we have
\[
I_\la(su_\ep)\le\frac{s^2}{2}\|\nabla u_\ep\|_2^2-\frac{s^q}{q}B(u_\ep)+\frac{4\la_0}{p}T^pA(w)\ep^\theta .
\]
By \eqref{eq:level}, the right-hand side is at most $\Theta-c_1\ep^\alpha+C\la_0\ep^\theta+o(\ep^\alpha)$. Moreover, for $s\in[T_0,T]$ we have $s^{q-2}\ge4q\mu$, so that
\[
\frac{s^2}{2}\|\nabla u_\ep\|_2^2-\frac{s^q}{q}B(u_\ep)\le s^2\mu-\frac{s^q}{2q}\le-\frac{s^q}{4q}.
\] Since $\theta>\alpha$, we may fix $\ep$ so small that
\begin{equation}\label{eq:path}
\max_{0\le s\le T}I_\la(su_\ep)\le\bar c:=\Theta-\frac{c_1}{2}\ep^\alpha,\quad I_\la(su_\ep)<0\ \text{for}\ s\in[T_0,T],\ \la\in(0,2\la_0].
\end{equation}

Following \cite{CZZ}, for $u\in H^1_0(\Om)$ and $v=\la^{1/(p-2)}u$ we have $I_\la(u)=\la^{-2/(p-2)}\Upsilon_m(v)$ and $I_\la'(u)=\la^{-1/(p-2)}\Upsilon_m'(v)$, where
\[
\Upsilon_m(v)=\frac12\|\nabla v\|_2^2+\frac1pA(v_+)-\frac mqB(v_+),\qquad m=m(\la):=\la^{-\frac{q-2}{p-2}} .
\]
We apply Jeanjean's theorem with $X=H^1_0(\Om)$, $\mathcal A(v)=\frac12\|\nabla v\|_2^2+\frac1pA(v_+)$, $\mathcal B(v)=\frac1qB(v_+)\ge0$, $J=[m(2\la_0),m(\la_0/2)]$, $v_1=0$ and $v_2=(2\la_0)^{1/(p-2)}T_0u_\ep$. For $m=m(\la)\in J$, we have $\la^{-1/(p-2)}v_2=s_\la u_\ep$ with $s_\la=(2\la_0/\la)^{1/(p-2)}T_0\in[T_0,T]$, so that \eqref{eq:path} gives $\Upsilon_m(v_2)=\la^{2/(p-2)}I_\la(s_\la u_\ep)<0=\Upsilon_m(0)$. On the other hand, by the Hardy--Sobolev inequality, $\Upsilon_m(v)\ge\frac12\|\nabla v\|_2^2-C\|\nabla v\|_2^q\ge\frac14r_1^2$ on $\{\|\nabla v\|_2=r_1\}$ for some $r_1>0$ independent of $m\in J$, and every path in $\mathcal G$ crosses this sphere. Hence $\tilde c_m\ge\frac14r_1^2>0$. Jeanjean's theorem yields, for almost every $m\in J$, a bounded Palais--Smale sequence $(v_k)$ of $\Upsilon_m$ at the level $\tilde c_m$. Then $u_k=\la^{-1/(p-2)}v_k$ is a bounded Palais--Smale sequence of $I_\la$ at the level $c_\la=\la^{-2/(p-2)}\tilde c_m>0$. Using the path $t\mapsto ts_\la u_\ep\la^{1/(p-2)}$ in $\mathcal G$ and \eqref{eq:path}, we obtain
\[
c_\la\le\max_{0\le s\le s_\la}I_\la(su_\ep)\le\bar c .
\]
Since $m\mapsto\la$ is a smooth decreasing bijection, the set of such $\la$ has full measure in $[\la_0/2,2\la_0]$.
\end{proof}

\begin{lemma}\label{lem:GK}
Assume $(H_\alpha)$. If $u_n\rightharpoonup0$ in $H^1_0(\Om)$, then $\limsup_nB(u_n)^{2/q}\le\mu^{-1}\limsup_n\|\nabla u_n\|_2^2$.
\end{lemma}

\begin{proof}
On $\Om\setminus B_\delta$ the weight is bounded and $q<2^*$, so $\int_{\Om\setminus B_\delta}|x|^{-s_2}|u_n|^q\to0$. In $B_{2\delta}$ we flatten the boundary by $\varphi$, which satisfies $|D\varphi-I|\le C\delta^\alpha$ and $|\varphi(x)|/|x|\in[1-C\delta^\alpha,1+C\delta^\alpha]$, and apply the Hardy--Sobolev inequality in $\R^N_+$ to a cut-off of $u_n$. This gives the claim up to a factor $1+C\delta^\alpha$, and we let $\delta\to0$, see also \cite{GK}.
\end{proof}

\begin{proof}[Proof of Theorem \ref{thm:main}(i)]
Fix $\la_0>0$. Let $(v_k)$ be a bounded Palais--Smale sequence of $I_\la$ at a level $c\in(0,\Theta)$ and $v_k\rightharpoonup v$. If $v=0$, then $A(v_k)\to0$ since $p<2^*(s_1)$, $\|\nabla v_k\|_2^2-B(v_k)\to0$ and $(\frac12-\frac1q)\|\nabla v_k\|_2^2\to c$. By Lemma \ref{lem:GK}, $L=\lim\|\nabla v_k\|_2^2$ satisfies $L^{2/q}\le\mu^{-1}L$, which gives $c\ge\Theta$, a contradiction. Hence $v\ge0$ is a nontrivial critical point, and $v>0$ by the strong maximum principle. Since
\[
I_\la(u)-\frac12\langle I_\la'(u),u\rangle=\Big(\frac12-\frac1q\Big)B(u_+)-\la\Big(\frac12-\frac1p\Big)A(u_+),
\]
$B$ is weakly lower semicontinuous and $A$ is compact, we also have $I_\la(v)\le c$.

By Lemma \ref{lem:ps}, there are $\la_n\to\la_0$, $\la_n\in[\la_0/2,2\la_0]$, and solutions $u_n$ of $(P_{\la_n})$ with $I_{\la_n}(u_n)\le\bar c<\Theta$. By the Nehari identity and Corollary \ref{cor:univ} with $\la_-=\la_0/2$,
\[
\Big(\frac12-\frac1q\Big)\|\nabla u_n\|_2^2=I_{\la_n}(u_n)+\la_n\Big(\frac1q-\frac1p\Big)A(u_n)\le\bar c+2\la_0C_A .
\]
Up to a subsequence, $u_n\rightharpoonup u_0$ in $H^1_0(\Om)$, and by Proposition \ref{prop:univ} and interior estimates, $u_n\to u_0$ in $C^1_{\rm loc}(\Om)$. Thus $u_0\ge0$ is a weak solution of $(P_{\la_0})$. If $u_0=0$, then $A(u_n)\to0$, and the Nehari identity together with the Hardy--Sobolev inequality on $\Om$ gives $\|\nabla u_n\|_2^2\ge\mu_{s_2}(\Om)^{q/(q-2)}>0$. Arguing as above, we obtain $I_{\la_n}(u_n)\ge\Theta-o(1)$, a contradiction. Therefore $u_0>0$ solves $(P_{\la_0})$.
\end{proof}

\section{Finiteness of \texorpdfstring{$\la^*$}{lambda*} and the behaviour at \texorpdfstring{$\la^*$}{lambda*}}\label{sec:finite}

Throughout this section we assume $(H_\alpha)$. Let 
\begin{align*}
	\mathcal L_\gamma=-\Delta-\gamma|x|^{-2}, \quad  \beta=\frac N2+\alpha \quad \text{and} \quad 0<\gamma\le\frac{a^2}{2}. 
\end{align*}
Note that $(\beta-1-\alpha)(N-1+\alpha-\beta)=a^2$.

\begin{lemma}\label{lem:barrier}
Let $r_0\le r_\Om$. There exist $M>0$ and $\Psi\ge0$ with $\Psi\in H^1(\Om\setminus B_r)$ for every $r>0$, $\mathcal L_\gamma\Psi\ge0$ weakly in $\Om$, $\Psi\le C$ in $\Om\setminus B_{r_0/4}$, and
\begin{gather*}
\Psi=\Psi_0:=(x_N+M|x|^{1+\alpha})|x|^{-\beta}\quad\text{in}\ \Om\cap B_{r_0/4},\\
d|x|^{-\beta}\le\Psi_0\le(1+M)|x|^{1-\beta}\quad\text{in}\ \Om\cap B_{r_0}.
\end{gather*}
\end{lemma}

\begin{proof}
For $M\ge K$ we have $x_N+M|x|^{1+\alpha}\ge x_N-h(x')\ge d$ in $\Om\cap B_{r_0}$. Since
\[
-\Delta(x_N|x|^{-\beta})=\beta(N-\beta)x_N|x|^{-\beta-2},\qquad-\Delta|x|^{1+\alpha-\beta}=a^2|x|^{\alpha-1-\beta},
\]
we obtain
\[
\mathcal L_\gamma\Psi_0=(\beta(N-\beta)-\gamma)x_N|x|^{-\beta-2}+M(a^2-\gamma)|x|^{\alpha-1-\beta}.
\] Both terms are nonnegative where $x_N\ge0$. Where $x_N<0$ we have $|x_N|\le K|x|^{1+\alpha}$, so that 
\begin{align*}
	 \mathcal L_\gamma\Psi_0\ge0 \quad \text{provided} \quad  M\ge2\beta(N-\beta)K/a^2.
\end{align*}
 Let $\chi\in C^\infty_c(B_{r_0})$ with $\chi=1$ on $B_{r_0/2}$, and $f=(\mathcal L_\gamma(\chi\Psi_0))_-$, which is bounded and supported in $B_{r_0}\setminus B_{r_0/2}$. Since $\gamma<a^2$, there is 
 \begin{align*}
 	\tilde w\in H^1_0(\Om\setminus\overline B_{r_0/4})\cap L^\infty, \quad \tilde w\ge0, \quad \text{with}~~\mathcal L_\gamma\tilde w=f.
 \end{align*}
   By the zero extension property in the proof of Lemma \ref{lem:KO}, its extension by zero satisfies $\mathcal L_\gamma\tilde w\ge f$ in $\Om$. Then $\Psi=\chi\Psi_0+\tilde w$ has the required properties.
\end{proof}

Fix $\eta\in(0,c_0)$ with $\eta^{q-2}\le\frac{a^2}{2}$, set $\gamma=\eta^{q-2}$, and define
\[
\rho(u)=\sup\{|x|:x\in\Om,\ \Phi_u(x)\ge\eta\},
\]
which is well defined by Lemma \ref{lem:nondeg}. By Lemma \ref{lem:KO}, $\Phi_u(x)\ge\eta$ implies $L_\la(x)\le C_\eta$. Thus 
\begin{equation}\label{eq:rho}
\la\rho(u)^\theta\le C_\eta .
\end{equation}

\begin{lemma}\label{lem:decay}
Let $\la\ge\la_-$ and $\rho=\rho(u)\le r_0/16$. Then $u\le C\rho^\alpha\Psi$ in $\Om\setminus B_{2\rho}$, and
\[
u\le C\rho^\alpha|x|^{1-\beta},\qquad|\nabla u|\le C\rho^\alpha|x|^{-\beta}\qquad\text{in}\ \Om\cap(B_{r_0/8}\setminus B_{4\rho}).
\]
\end{lemma}

\begin{proof}
For $|x|>\rho$ we have $\Phi_u<\eta$ and  $-\Delta u\le\gamma|x|^{-2}u$. On $\Om\cap\partial B_{2\rho}$,  Lemma \ref{lem:scaled} gives $|\nabla u|\le C\rho^{-a-1}$, and integrating along normal segments we obtain $u\le Cd\rho^{-a-1}\le C'\rho^\alpha\Psi$, since $\beta-a-1=\alpha$. Let $w=u-C'\rho^\alpha\Psi$ and $v_\ep=(w-\ep)_+1_{\{|x|>2\rho\}}$. Since $0\le v_\ep\le u$ and $v_\ep$ vanishes near $\partial B_{2\rho}$, we have $v_\ep\in H^1_0(\Om\setminus\overline B_{2\rho})$. Testing $\mathcal L_\gamma w\le0$ with $v_\ep$ and using the Hardy inequality, we get
\[
\Big(1-\frac{\gamma}{a^2}\Big)\int|\nabla v_\ep|^2\le\gamma\ep\int_{\Om\setminus B_{2\rho}}|x|^{-2}u\to0 ,
\]
so that $w\le0$. The other bounds follow from $\Psi_0\le(1+M)|x|^{1-\beta}$ and Lemma \ref{lem:scaled}.
\end{proof}

\begin{lemma}\label{lem:Alow}
Let $\la_n\ge\la_-$, $\rho_n=\rho(u_n)\to0$ and fix $r>0$. Then $\liminf\limits_n\rho_n^{-\theta}A(u_n,B_r)>0$.
\end{lemma}

\begin{proof}
Let $w_n(y)=\rho_n^au_n(\rho_ny)$. Then $w_n$ solves the equation with $\la$ replaced by $\la_n\rho_n^\theta\le C_\eta$, $\Phi_{w_n}\le C_0$, and the domains $\rho_n^{-1}\Om$ converge to $\R^N_+$. By Lemma \ref{lem:scaled}, $w_n\to W$ in $C^1_{\rm loc}(\overline{\R^N_+}\setminus\{0\})$. Pick $x_n$ with $\rho_n/2\le|x_n|\le\rho_n$ and $\Phi_{u_n}(x_n)\ge\eta$. Since $w_n\le Cd$ on $\{\frac14\le|y|\le2\}$, the points $x_n/\rho_n$ stay away from $\partial\R^N_+$, and  $W\not\equiv0$ on $\{\frac14\le|y|\le2\}$. The claim follows from $A(u_n,B_r)=\rho_n^\theta A(w_n,B_{r/\rho_n})$.
\end{proof}

\begin{proposition}\label{prop:key}
Assume $\theta\le\alpha$ and $\la_->0$. There exist $\rho_0,C>0$ such that every solution of \eqref{eq:P} with $\la\ge\la_-$ and $\rho(u)\le\rho_0$ satisfies $\la\le C\rho(u)^{\alpha-\theta}$.
\end{proposition}

\begin{proof}
Let $r=r_0/8$. By Lemma \ref{lem:Alow} and a contradiction argument, there exist $\rho_0\le r_0/16$ and $c_1>0$ such that $A(u,B_r)\ge c_1\rho^\theta$ whenever $\rho=\rho(u)\le\rho_0$. By \eqref{eq:poh},
\[
c_pc_1\la\rho^\theta\le|P(r,u)|+\frac12\int_{\partial\Om\cap B_r}|x\cdot\nu||\partial_\nu u|^2 .
\]
On $\Om\cap\partial B_r$, Lemma \ref{lem:decay} gives $u,|\nabla u|\le C\rho^\alpha$, and with \eqref{eq:rho} we get
\[
|P(r,u)|\le C(\rho^{2\alpha}+\la\rho^{p\alpha})\le C\rho^\alpha,
\]
since $(p-1)\alpha>\theta$. By \eqref{eq:xnu}, Corollary \ref{cor:univ} and Lemma \ref{lem:decay},
\begin{align*}
\int_{\partial\Om\cap B_r}|x\cdot\nu||\partial_\nu u|^2&\le C\rho^\alpha+C\rho^{2\alpha}\int_{4\rho}^rs^{1+\alpha-2\beta}s^{N-2}ds\\
&=C\rho^\alpha+C\rho^{2\alpha}\int_{4\rho}^rs^{-1-\alpha}ds\le C\rho^\alpha .
\end{align*}
Hence $\la\rho^\theta\le C\rho^\alpha$.
\end{proof}

\begin{proof}[Proof of Theorem \ref{thm:main}(ii)]
If $\la_n\to\infty$ along solutions $u_n$, then \eqref{eq:rho} gives $\rho(u_n)\to0$, and Proposition \ref{prop:key} gives $\la_n\le C\rho(u_n)^{\alpha-\theta}\le C$, a contradiction.
\end{proof}

We now turn to $\la^*$.

\begin{proposition}\label{prop:star}
Assume $\theta<\alpha$ and $0<\la^*<\infty$. Let $\la_n\to\la^*$ and let $u_n$ solve $(P_{\la_n})$. Then the following hold.
\begin{enumerate}[label=(\roman*)]
\item Up to a subsequence, $u_n\to u_0$ in $C^1_{\rm loc}(\Om)$, where $u_0\not\equiv0$ solves the equation in $\mathcal D'(\Om)$ and $\Phi_{u_0}\le C_0$.
\item Either 
\begin{align*}
&\limsup_{x\to0}\Phi_{u_0}(x)=0 \quad \text{and} \quad u_0\in H^1_0(\Om) \quad \text{solves} \quad  (P_{\la^*}) \quad \text{or} \\
&  \limsup_{x\to0}\Phi_{u_0}(x)\ge c_0 \quad \text{and} \quad u_0\notin H^1_0(\Om).
\end{align*}
\item $\la^*$ is attained if and only if there are solutions $u_n$ of $(P_{\la_n})$, $\la_n\to\la^*$, with $\sup\limits_n\|\nabla u_n\|_2<\infty$.
\end{enumerate}
\end{proposition}

\begin{proof}
(i) By Proposition \ref{prop:key} with $\la_-=\la^*/2$ and $\theta<\alpha$, $\rho_n:=\rho(u_n)\ge\rho_*>0$. Pick $x_n$ with $|x_n|\ge\rho_*$ and $u_n(x_n)\ge\eta(\operatorname{diam}\Om)^{-a}$. On $\{|x|>\rho_*/4\}$ we have $u_n\le C$ and $|V_n|\le C$, where $-\Delta u_n=V_nu_n$. Testing with $u_n\zeta^2$, where $\zeta=0$ in $B_{\rho_*/4}$ and $\zeta=1$ outside $B_{\rho_*/2}$, we see that $(u_n\zeta)$ is bounded in $H^1_0(\Om)$, so that $u_n\zeta\to u_0\zeta$ in $L^2(\Om)$. The zero extension of $u_n$ satisfies $\Delta\tilde u_n\ge-C\tilde u_n$ in $\{|x|>\rho_*/2\}$, and the local maximum principle \cite[Theorem 8.17]{GT} gives $u_n(x_n)\le C\|u_n\|_{L^2(B(x_n,\rho_*/4))}$. Thus $u_0\not\equiv0$. No regularity of $\partial\Om$ away from $0$ is used.

(ii) Suppose $\Phi_{u_0}^{q-2}\le\gamma'<a^2$ in $B_r$. Let $\zeta_\ep=0$ in $B_\ep$, $\zeta_\ep=1$ outside $B_{2\ep}$, $|\nabla\zeta_\ep|\le2/\ep$. Testing the equation of $u_0$ with $u_0\zeta_\ep^2$, we obtain
\[
\Big(1-\frac{\gamma'}{a^2}\Big)\int|\nabla(u_0\zeta_\ep)|^2\le C\int u_0^2+\frac{4}{\ep^2}\int_{B_{2\ep}\setminus B_\ep}C_0^2|x|^{-2a}\le C,
\]
since $N-2a=2$. Hence $u_0\in H^1_0(\Om)$ is a weak solution of $(P_{\la^*})$, and Lemma \ref{lem:reg}(ii) gives $\Phi_{u_0}\le C|x|^{a+1}\to0$. Conversely, if $u_0\in H^1_0(\Om)$, then $u_0$ solves $(P_{\la^*})$ and Lemma \ref{lem:reg}(ii) gives $\Phi_{u_0}\to0$. Since $c_0^{q-2}=a^2$, this proves the dichotomy.

(iii) If $\sup\limits_n\|\nabla u_n\|_2<\infty$, then $u_n\rightharpoonup u_0$ in $H^1_0(\Om)$, and $u_0\not\equiv0$ solves $(P_{\la^*})$ by (i). The converse is obvious.
\end{proof}

\begin{remark}\label{rem:open}
	By Proposition \ref{prop:star}, if $\theta<\alpha$, then $\la^*$ fails to be attained if and only if, for every sequence $\la_n\uparrow\la^*$ and every choice of solutions $u_n$ of $(P_{\la_n})$, the limit $u_0$ has a nonremovable singularity at $0$, that is, $\limsup\limits_{x\to0}\Phi_{u_0}(x)\ge c_0$ and $u_0\notin H^1_0(\Om)$. 
	One possible mechanism is an infinite sequence of concentration scales $r_k\to0$. A formal Pohozaev balance between consecutive scales gives
	\[
	c\,\la^*r_k^{\frac N2}\le r_{k+1}^{\frac N2-\theta}\le C\,\la^*r_k^{\frac N2}\qquad\text{for all}\ k,
	\]
	where the constants $0<c\le C$ do not depend on $k$. This is compatible with \eqref{eq:poh}.
	 We do not know whether this scenario, or any other nonremovable singularity, occurs. The question is thus whether the perturbed Hardy--Sobolev equation
	\[
	-\Delta u=|x|^{-s_2}u^{q-1}-\la^*|x|^{-s_1}u^{p-1}
	\]
	admits solutions with a nonremovable isolated boundary singularity at $0$ that arise as such limits. For the unweighted critical equation $-\Delta u=u^{\frac{N+2}{N-2}}$, isolated boundary singularities were analysed by Xiong \cite{X}, and Wang and Xiong \cite{WX} showed that their asymptotic behaviour need not be stationary. The weighted case $s_2>0$, which is the relevant one here, will be studied elsewhere.
	
	In the doubly critical case $p=2^*(s_1)$, that is, $\theta=0$, the finiteness of $\la^*$ is proved in \cite{TT1}. In our notation it follows from the identities
	\[
	|x|^{-s_2}u^q=|x|^{-2}u^2\Phi_u^{q-2},\qquad|x|^{-s_1}u^p=|x|^{-2}u^2\Phi_u^{p-2},
	\]
	the Nehari identity and the Hardy inequality, which give $\la^*\le\la_c$, where $\la_c$ is the unique $\la$ with $\max\limits_{t>0}(t^{q-2}-\la t^{p-2})=a^2$. In this case the equation is invariant under $u\mapsto\mu^au(\mu\,\cdot)$ with $\la$ fixed, Step 2 in the proof of Proposition \ref{prop:univ} has no analogue, and the attainability of $\la^*$ remains open.
\end{remark}

\end{document}